# In Memoriam
# Vadim Anatol'evich Yankov (1935-2024)

ALEX CITKIN
*Metropolitan Telecommunications, New York, USA*

IOANNIS M. VANDOULAKIS
*Open University of Cyprus, Cyprus & Logica Universalis Association*

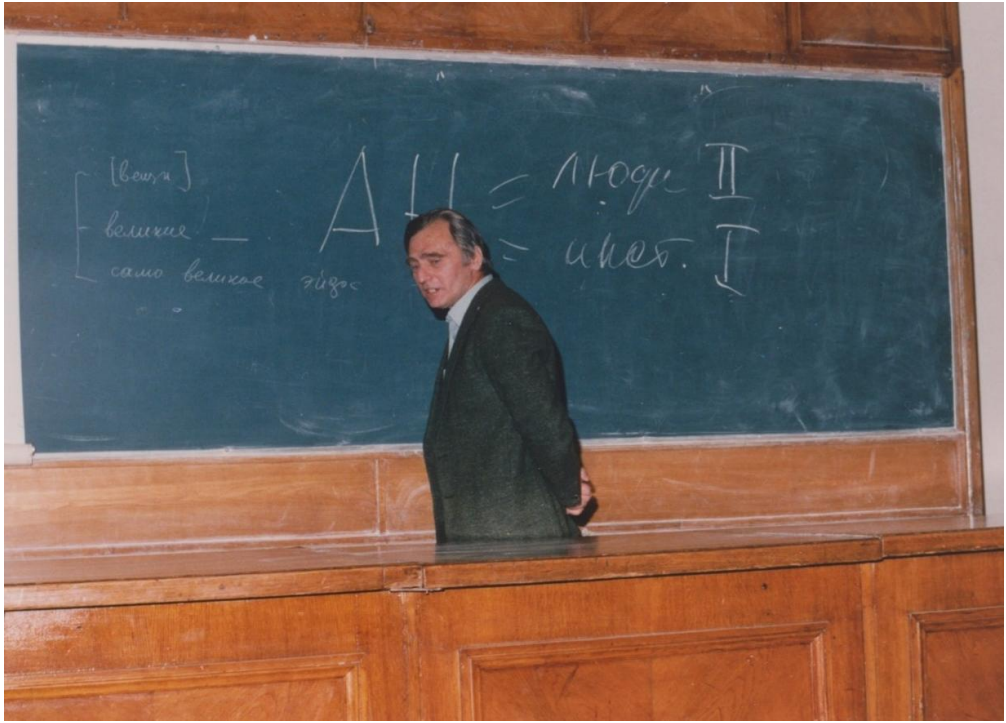

*Figure 1. Yankov, Moscow Lomonosov University, 4 October 1991.*
*Photo by I. Vandoulakis, first published in [Vandoulakis 2015].*

On March 11, 2024, the Russian mathematical logician, philosopher, historian of logic and the foundations of mathematics and political activist who was prosecuted in the former USSR, Vadim Anatol'evich Yankov (or Jankov), passed away in Dolgoprudny, a town in Moscow Oblast, Russia, at the age of 89. His health was seriously affected by at least four strokes that left him in bed for a long time. His death signals something like the end of a chapter in the history of Markov's school of constructive mathematics.

He was born in Taganrog, a city located on the north shore of Taganrog Bay in the Sea of Azov, into an educated family. His father, Anatoly, was an engineer from Russia, and his mother was a schoolteacher of Hebrew origin but professed Orthodoxy. He spent his childhood in Taganrog until 1941 when the Second World War compelled the family to evacuate to Sverdlovsk. After the end of the war, in 1943, the family returned to Taganrog, which Vadim Anatol'evich left in 1952 for Moscow, where he enrolled in the Department of Philosophy of the Moscow M.V. Lomonosov State University. However, he soon realised that the ideologisation in the humanities during the Soviet time did not leave enough room for free thinking. Thus, one year later, he decided to change to the Department of Mechanics and Mathematics. However, he soon faced troubles there, too, because of his critical political mind. In 1956, he was expelled from the University for his sharp criticism of the Komsomol, his participation in a complaint against the conditions at the University students' restaurant and publication of an independent students' newspaper. Later, he was accepted into the University's distance learning program and obtained his diploma in 1959.

One year before his graduation, he was employed in the Programming Department at the Steklov Institute of Mathematics (later, the Programming Department of the Institute of Mathematics of the Siberian Branch of the USSR Academy of Sciences). He joined the Research Group, which was developing one of the first programming languages, the ALPHA, an extension of ALGOL.

He was married in 1960 to Natalya Sarmakesheva, of Armenian origin. They had two sons, Kirill and Ilya, and a daughter, Anastasiya.

After his graduation in 1959, he was accepted by Andrey Andreevich Markov (1903-1979), the founder of the Russian School of Constructive Mathematics and Logic, as a post-graduate student at the Faculty of Mathematics of Moscow State University. Markov had been elected the Head of the Department of Mathematical Logic at Moscow University this year and was also a member of the Mathematical Institute of the Academy of Sciences (1939-1972). Under Markov's supervision, Yankov completed his PhD Thesis "Finite Implicative Structures and Realizability of Formulas of Propositional Logic" and was awarded his PhD in 1964.

During his PhD studies, Yankov started his research in Propositional Logic. Within a few years, he published several papers which revolutionised the theory of Propositional Logics, specifically the Intermediate Logics. His most important observation was that the propositional language of intermediate logics is rich enough to describe some properties of finite algebraic models, the so-called *pseudo-Boolean (Heyting) algebras*: with any finite irreducible Heyting algebra **A,** one can associate a propositional formula $\chi(\mathbf{A})$ and $\chi(\mathbf{A})$ is refuted in Heyting algebra **B** if and only if algebra **A** can be embedded in a homomorphic image of **B**. It turned out that characteristic formulas, and only them, are conjunctively irreducible formulas of the Intuitionistic Propositional Logic (IPC). Moreover, using the characteristic formulas, Yankov proved that there is a continuum of many intermediate logics, and there exist intermediate logics without the finite model property. He paid special attention to the intermediate logic defined by the weak law of excluded middle: $\neg p \lor \neg\neg p$, and he proved that all intermediate logics having the same positive fragments as IPC lay between IPC and this logic, which nowadays is known as the *Yankov Logic*. His papers of this period are the most significant and most cited papers in the field of propositional and algebraic logics [Citkin 2022].

In 1963, Yankov started teaching at the Moscow Institute of Physics and Technology, where he was employed as an Assistant Lecturer. However, he was dismissed in 1968, when he co-signed the famous "Letter of the 99 Soviet mathematicians" addressed to the Ministry of Health and the General Procurator of Moscow asking for the release of Aleksandr Esenin-Vol'pin [Fuchs 2007, 221].

Between 1968 and 1974, Yankov worked as a Senior Lecturer at the Moscow Aviation Institute, where he was again fired for his frequent dissenting discussions about the Soviet intervention in Czechoslovakia in 1968. Thus, in 1974, he found a job at the Enterprise Resource Planning Department of the Moscow Institute for Urban Economics, where he worked until 1982.

From November 1981 until January 1982, he published some critical papers abroad in the dissident journal *Kontinent,* founded in 1974 by writer Vladimir Maximov, which was printed in Paris and focused on the politics of the Soviet Union. In issue 18, he published the article "On the Possible Meaning of the Russian Democratic Movement", and in 1981-82, the article "Letter to Russian Workers about the Polish Events," expressing his support for the broad anti-authoritarian social movement of the Polish trade union "Solidarity" (*Solidarność*). Thus, on August 9, 1982, when Yankov left his apartment to go to the office, he was arrested. In January 1983, he was sentenced by the Moscow City Court to four years in prison and three years of exile for anti-Soviet propaganda. He served his term in the Gulag labour camp for political prisoners, called Dubravny Camp in Barashevo, Mordovia, and exile in Buryatia in south-central Siberia.

While in prison, he did not lose his positive spirit. He gathered up the courage to start studying Classic Greek by comparing Thucydides' works in the original and their Russian translation and guided by Sergey Ivanovich Sobolevsky's (1864–1963) textbook on Classic Greek language, first published in 1948, with which a whole generation of classicists grew up in the Soviet Union. According to reminiscences from his fellow imprisoned inmate Levan Berdzenishvili, a Georgian politician and academic who, after his release, served as the director of the National Parliamentary Library of Georgia from 1998 until 2004, he liked to organise "Socratic dialogues" among them, where he performed the role of Socrates or an experienced sophist [Berdzenishvili [no date]]. Yankov appeared in 2019 as the main character in Berdzenishvili's novel *Sacred Darkness*, describing their joint time in prison. [Berdzenishvili 2019]

During the *Perestroika*, in January 1987, he was released and rehabilitated on October 30, 1991, from all charges against him. After his release, he worked in the Institute of Thermal Metallurgical Units and Technologies "STALPROEKT" until 1991.

Since 1991, he has been an Associate Professor at the Department of Mathematics, Logic and Intellectual Systems, Faculty of Theoretical and Applied Linguistics and the Department of Mathematical and Logical Foundations of Humanitarian Knowledge, Institute of Linguistics of the Russian State University for the Humanities in Moscow that was founded in 1991 by merging the Moscow Urban University of the People (est. 1908) and the Moscow State University for History and Archives (est. 1930) and envisaged as a new type of university focused on reviving the liberal arts in Russia that had suffered from ideologisation during the Soviet period.

This was the second fruitful period of his intellectual life. Although he initially taught standard mathematical courses, such as calculus, algebra, probability, and computer science, his research and educational interests gradually shifted to philosophy, the history of philosophy, and the history of mathematics. Thus, he started lecturing regular courses on philosophy and the history of philosophy at the Russian State University of the Humanities, delivered a series of lectures in the Seminar of Philosophy of Mathematics of the Moscow M.V. Lomonosov State University and published papers in philosophy and the history of mathematics.

Yankov's philosophical concerns started being shaped while he was imprisoned. His first, possibly philosophical publication was printed abroad in issue 43 (1985) of the journal *Kontinent*, entitled "Ethical-philosophical treatise," where he outlines his philosophical conception of existential history. A publication on the same theme in Russia was made possible only ten years later in the journal *Voprosy Filosofii* [Yankov 1998]. The reception of his ideas in this article was extremely critical, partly because of his chaotic style of exposition. Nevertheless, his idea of linking different cultures into a single semantic whole and tracking their development in time, taking the existential idea as the central axis, was later appreciated as innovative and promising [Denisova 2020, 272].

His philosophical inquiry culminated in the 2011 publication of his capital monograph, *An Interpretation of Early Greek Philosophy* [Yankov 2011]. In it, he displays a panoramic view of the early Greek picture of the world, tracing the existential line as an accidental and optional element of a generally "faithful" interpretation [Denisova 2020, 274]. Yankov's most significant contribution is that he highlighted the idea of anthropological problematics in the early Greek philosophy, which has been stated before him as a vague hypothesis of some deviant interpretive trend, into a research problem. [Denisova 2020, 292].

At the same time, his book is a profound study of the beginnings of mathematics, the appearance of philosophy, and rational thinking generally. Yankov focuses on the appearance of proof in ancient Greece, a distinctive feature of the European scientific tradition, a topic he studied for almost a decade [Yankov 1997, 2000, 2001, 2003]. It is noteworthy that Yankov examines the appearance of proof in Greek mathematics in relationship to the development of ontology in early Greek philosophy as phenomena that concern the rise of rational thinking. His profound insights provide a new, consistent narrative of the historical development of logical thinking in mathematics and philosophy and their interrelationship. [Vandoulakis 2022]

Yankov also possesses an important place in the history of foundational studies in the Soviet Union because he participated in the foundational debates of that time. As a disciple of A.A. Markov, he witnessed the development of Markov's constructive school and, particularly, contributed to preserving Markov's views about Brouwer's intuitionism and Hilbert's formalism. He achieved a significant historical record of Markov's views, up to 1965, by including Markov as a peculiar interlocutor of Brouwer and Hilbert in the Russian translation of Heyting's *Intuitionism* [Vandoulakis 2015, 2022]. Markov's dialogue with Heyting's fictional persons initiated a real dialogue between Markov and the intuitionists when Heyting responded to it in the third edition of his book in 1971. Moreover, it offers the necessary historical material to reconstruct the principles underlying Markov's "philosophy of constructivism."

The volume [Citkin, Vandoulakis (Eds) 2022] dedicated to him in the series Outstanding Contributions to Logic (OCTR, volume 24) and the Logica Universalis Webinar [Yankov, Citkin, Vandoulakis and Denisova 2024] organised to his memory are a minimal appreciation to a mathematician and scholar who deserves our respect and admiration. A complete list of Yankov's contributions can be found in [Citkin, Vandoulakis 2022, 3-4].

## Acknowledgement

The authors wish to express their appreciation to Kirill V. Yankov, son of Prof. V.A. Yankov, for the information about his life and the chronology of some of his life's events.